\input ./style/arxiv-general.cfg
\documentclass[MSNbibl,number,citesort,seceqn,dvips]{arxbj}
\makeatletter
   \@ifpackageloaded{graphicx}{}{\usepackage{graphicx}}
\makeatother
\usepackage{upgreek}

\aid{0}
\volume{21}
\issue{3}
\pubyear{2015}
\firstpage{1341}
\lastpage{1360}
\doi{10.3150/13-BEJ588} % kopijuoti is 'New paper accepted'
\makeatletter
\newcommand{\angler}{\rangle}
\newcommand{\anglel}{\langle}
\newcommand{\rrvert}{\vert}
\newcommand{\llvert}{\vert}
\newtheorem{theorem}{Theorem}[section]
\newtheorem{lemma}{Lemma}[section]
\newtheorem{proposition}{Proposition}[section]
\newtheorem{corollary}{Corollary}[section]
\newcommand{\eqref}[1]{(\ref{#1})}
\newremark{remark}{Remark}[section]
\makeatother

\begin{document}
\begin{frontmatter}

\title{Mimicking self-similar processes}
\runtitle{Mimicking self-similar processes}

\begin{aug}
%%%% inicialai - be tarpu
\author[A]{\inits{J.Y.}\fnms{Jie Yen} \snm{Fan}\thanksref{e1}\ead[label=e1,mark]{jieyen.fan@monash.edu}},
\author[A]{\inits{K.}\fnms{Kais} \snm{Hamza}\thanksref{e2}\corref{}\ead[label=e2,mark]{kais.hamza@monash.edu}}
\and
\author[A]{\inits{F.}\fnms{Fima} \snm{Klebaner}\thanksref{e3}\ead[label=e3,mark]{fima.klebaner@monash.edu}}
%%\runauthor{} %% auto
%\dedicated{}
\address[A]{School of Mathematical Sciences, Monash University, Clayton,
VIC 3800, Australia. \\\printead{e1,e2,e3}}
%\address[]{}
\end{aug}

% HISTORY:
\received{\smonth{9} \syear{2012}}
\revised{\smonth{8} \syear{2013}}

% ABSTRACT
%
\begin{abstract}
We construct a family of self-similar Markov martingales
with given marginal distributions. This construction uses the
self-similarity and Markov property of a reference process to
produce a family of Markov processes that possess the same marginal
distributions as the original process. The resulting processes are
also self-similar with the same exponent as the original process. They
can be chosen to be martingales under certain conditions. In this paper,
we present two approaches to this construction, the transition-randomising
approach and the time-change approach. We then compute the infinitesimal
generators and obtain some path properties of the resulting processes.
We also give some examples, including continuous Gaussian martingales
as a generalization of Brownian motion, martingales of the squared
Bessel process, stable L\'evy processes as well as an example of an
artificial process having the marginals of $t^\kappa V$ for some
symmetric random variable $V$. At the end, we see how we can mimic
certain Brownian martingales which are non-Markovian.
\end{abstract}

% KEYWORDS
% visi is mazosios raides ir pagal abecele
%

\end{frontmatter}

%s1 #&#
\section{Introduction}

Constructing martingales with given marginal distributions has been an
active area of research over the last decade (e.g. \cite{Alb08,BakDonYor11,HamKle07,HirProRoyYor11,MadYor02,Ole08}).
(Here and in the entire paper, marginal distributions (also marginals)
refer to the 1-dimensional distributions.)
A condition for the existence of such martingales is given by Kellerer
\cite{Kel72} (see Hirsch and Roynette \cite{HirRoy12} for a new and
improved proof).

Three constructions of Markov martingales with pre-specified marginal
distributions were given by Madan and Yor \cite{MadYor02}, namely the
Skorokhod embedding method, the time-changed Brownian motion and the
continuous martingale approach pioneered by Dupire \cite{Dup94}.
Recently, Hirsch \textit{et al}. \cite{HirProRoyYor11} gave six different
methods for constructing martingales whose marginal distributions match
those of a given family of probability measures. They also tackle the
tedious task of finding sufficient conditions to ensure that the chosen
family is indeed increasing in the convex order, or as they coined it,
a peacock.

In this paper, we deal with a different, albeit related, scenario. We
do not start with a family of probability distributions, rather we
start with a given martingale (the existence of which is assumed) and
produce a large family (as opposed to just a handful) of new
martingales having the same marginal distributions as the original process.
We say that these martingales ``mimic'' the original process.
%Note that, since we focus on mimicking existing processes, the
%question of the existence of a martingale with given marginals (those
%of the original martingale) is entirely redundant.

This same task was undertaken in \cite{HamKle07} for the Brownian
motion. It gave rise to the papers \cite{Alb08,BakDonYor11} and
\cite{Ole08} (who coined the term Faked Brownian Motion).
Albin \cite{Alb08} and Oleszkiewicz \cite{Ole08} answered the question
of the existence of a continuous martingale with Brownian marginals.
However, their constructions yield non-Markov processes.
Baker \textit{et al}. \cite{BakDonYor11} then generalised Albin's construction
and produced a sequence of (non-Markov) martingales with Brownian marginals.
In this paper, we extend the construction of \cite{HamKle07} to a much
larger class of processes, namely self-similar Markov martingales.

Before formulating a solution to this problem we give a brief account
on the origin and relevance of the mimicking question to finance, and
more specifically to option pricing; that is the pricing of a contract
that gives the holder the right to buy (or sell) the instrument (a
stock) at a future time $T$ for a specified price $K$. The theoretical
valuation of an option is performed in such a way as not to allow
arbitrage opportunities -- arbitrage occurs when riskless trading
results in profit.
The first fundamental theorem (e.g., \cite{Shr04}, page 231) states that
the absence of arbitrage in a market with stock price $S_t$, $0\le t\le
T$, is essentially equivalent to the existence of an equivalent
probability measure under which the stock price is a martingale. (Here
without loss of generality, we let the riskless interest rate be zero.)
The second fundamental theorem (e.g., \cite{Shr04}, page 232) implies
that the arbitrage-free price of an option is given by $\mathbb
{E}[(S_T-K)^+]$, where the expectation is taken under the equivalent
martingale measure. In the classical model of Black and Scholes the
stock price $S_t$ is given, under the martingale measure, by the
exponential Brownian motion $S_t=S_0\exp(\sigma B_t-t\sigma^2/2)$. The
parameter $\sigma$ is known as the volatility, and is assumed to be a
constant. The resulting expectation produces the well-known
Black--Scholes formula for option prices. However, empirical evidence
shows that in order to match the Black--Scholes formula to market prices
of options one needs to vary $\sigma$. As a consequence it is natural
to ask whether there exist alternative models of stock prices that
result in a prescribed option pricing formula, such as Black--Scholes.
Finally, it is easy to see that the collective knowledge of $\{\mathbb
{E}[(S_T-K)^+],K\geq0\}$ determines the distribution of $S_T$ or the
marginal distribution, for example \cite{HamKle06}. Therefore, if one
wants to keep the option prices given by the original formula but
without the limitations of the original process (such as constant
volatility) one has to look for martingales (to have the model
arbitrage free) with given marginals (to keep the same option prices).
This question received much attention in the last ten years, see the
pioneering work of Madan and Yor \cite{MadYor02}.

Throughout this paper, we assume that all processes are c\`adl\`ag and
progressively  measurable.
We will use the notation $\stackrel{d}{=}$ to mean equal in distribution for random variables or equal in finite-dimensional distributions for processes, and this will
be clear in the context.
For a given random measure $M(\mathrm{d}x)$, the measure
$M(c\,\mathrm{d}x)$ for $c>0$ is
defined by $\int g(x)M(c\,\mathrm{d}x) = \int g(x/c)M(\mathrm{d}x)$.
We will also write $\mathbb{E}[M(\mathrm{d}x)]$ to mean the
measure defined by
$\int g(x)\* \mathbb{E}[M(\mathrm{d}x)] = \mathbb{E}[\int g(x) M(\mathrm{d}x)]$ for any
positive function $g$.

We start with a martingale $Z$ which is also a Markov process, that is,
for any bounded measurable function $g$,
$\mathbb{E}[g(Z_t)|\mathcal{F}_s] = \mathbb{E}[g(Z_t)|Z_s]$ for $s\le t$,
where $(\mathcal{F}_t)_{t\ge0}$ denotes the natural filtration of
$(Z_t)_{t\ge0}$.
We aim to construct new processes that will have the same marginal
distributions as $Z$ while retaining the martingale and Markov properties.
We assume further that $Z$ is self-similar, that is, there exists a
(strictly) positive function $q(c)$ such that
$(Z_{ct})_{t\ge0} \stackrel{d}{=} (q(c)Z_t)_{t\ge0}$ for all $c>0$.
For $Z$ to be non-trivial and stochastically continuous at $t=0$, we
must have (see, e.g., \cite{EmbMae02})
that $q(c)=c^{\kappa}$ for some $\kappa>0$, i.e.,
\[
(Z_{ct})_{t\ge0} \stackrel{d} {=} \bigl(c^\kappa
Z_t\bigr)_{t\ge0}\quad\quad \forall c>0
\]
(which implies that $Z_0=0$).
We then say that $Z$ is self-similar with exponent $\kappa$, or simply,
$\kappa$-self-similar.

Denote the transition function of $Z$ by $P(s,t,x,\mathrm{d}y) := \mathbb
{P}(Z_t\in \mathrm{d}y| Z_s=x)$ for $s\le t$.
Then the self-similarity of $Z$ translates to the following scaling
property on $P$:
\[
P\bigl(cs,ct,c^\kappa x,c^\kappa\,\mathrm{d}y\bigr) = P(s,t,x,\mathrm{d}y)\quad\quad \forall
c>0, s\le t \mbox { and } x.
\]
If $L$ is the infinitesimal generator associated with $P$, this scaling
property is equivalent to
%
%e1.1 #&#
\begin{equation}
\label{E:SSGen} c\pi_{c^\kappa}L_{cs}\pi_{c^{-\kappa}} =
L_s,
\end{equation}
where $\pi_c$ is the operator defined by $\pi_{c}f(x) = f(cx)$.
From this, we see that
%
%e1.2 #&#
\begin{equation}
\label{E:SSGen1} L_s = s^{-1}\pi_{s^{-\kappa}}L_1
\pi_{s^{\kappa}}.
\end{equation}

In the following, we present a mimicking scheme to the process $Z$ by
randomising the transition functions.
We will see that this is equivalent to time-changing the process
together with an appropriate scaling.
We then obtain some properties of the resulting processes and give some
examples.
%At the end, we try to extend this construction to a more general class
%of processes.

%s2 #&#
\section{Mimicking scheme} \label{S:Scheme}

Let $Z$ be a Markov process with transition function $P$, which is
self-similar with exponent \mbox{$\kappa>0$}.
Note that if $Z$ is a martingale, then it has a c\`adl\`ag version;
if $Z$ is not a martingale, there are conditions for a Markov process
to be c\`adl\`ag.
In this section, we construct new Markov processes that possess the
same marginal distributions as $Z$ and show that the resulting
processes are martingales under certain conditions.

%le2.1 #&#
\begin{lemma} \label{LemmaP}
Let $0<s\le t\le u$. For any $a,b \in[0,1]$ and a measurable set $B$,
\[
\int P(0,s,0,\mathrm{d}x) P\bigl(at,t,(at/s)^{\kappa}x,B\bigr) = P(0,t,0,B)
\]
and
\[
\int P\bigl(at,t,(at/s)^{\kappa}x,\mathrm{d}y\bigr) P\bigl(bu,u,(bu/t)^{\kappa}y,B
\bigr) = P\bigl(abu,u,(abu/s)^{\kappa}x,B\bigr).
\]
\end{lemma}

\begin{pf}
Here we prove only the second equality, the first one is proved similarly.
Suppose that $a,b \in(0,1]$.
By the scaling property, we have
\begin{eqnarray*}
&&\int P\bigl(at,t,(at/s)^{\kappa}x,\mathrm{d}y\bigr) P\bigl(bu,u,(bu/t)^{\kappa}y,B
\bigr)
\\
&&\quad = \int P\bigl((bu/t)at,(bu/t)t,(bu/t)^{\kappa}(at/s)^{\kappa
}x,(bu/t)^{\kappa}\,\mathrm{d}y
\bigr) P\bigl(bu,u,(bu/t)^{\kappa}y,B\bigr)
\\
&&\quad = \int P\bigl(abu,bu,(abu/s)^{\kappa}x,\mathrm{d}y\bigr) P(bu,u,y,B) = P
\bigl(abu,u,(abu/s)^{\kappa}x,B\bigr).
\end{eqnarray*}
Notice that when $a$ or $b$ is $0$, the two equalities in the lemma
trivially hold true.
\end{pf}

%pr2.1 #&#
\begin{proposition}
Let $(G_{s,t})_{s\le t}$ be a family of probability distributions on
the set $(0,1]$, where $G_{s,s} = \delta_1$, the Dirac measure at $1$.
Suppose that for any bounded measurable function $h$ and $s\le t\le u$,
we have
%
%e2.1 #&#
\begin{equation}
\label{E:PropG} \int\int h(ab) G_{s,t}(\mathrm{d}a)G_{t,u}(\mathrm{d}b) = \int
h(r) G_{s,u}(\mathrm{d}r).
\end{equation}
Then $\widetilde{P}$ defined as follows is a transition function,
\begin{eqnarray*}
\widetilde{P}(0,t,0,\mathrm{d}y) &=& P(0,t,0,\mathrm{d}y),
\\
\widetilde{P}(s,t,x,\mathrm{d}y) &=& \int P\bigl(rt,t,(t/s)^{\kappa}r^{\kappa}x,\mathrm{d}y
\bigr) G_{s,t}(\mathrm{d}r),\quad\quad s\le t.
\end{eqnarray*}
\end{proposition}

\begin{pf}
Clearly, for each $(s,t,x)$, $\widetilde{P}(s,t,x,\mathrm{d}y)$ is a probability
measure and for each $(s,t,B)$,
$\widetilde{P}(s,t,x,B)$ is measurable in $x$. Note also that
$\widetilde{P}(s,s,x,B) = \delta_x(B)$.
Next, using Lemma~\ref{LemmaP}, we obtain, for $0<s\le t\le u$,
$\int\widetilde{P}(0,t,0,\mathrm{d}y) \widetilde{P}(t,u,y,B) = \widetilde{P}(0,u,0,B)$
and
\begin{eqnarray*}
\int\widetilde{P}(s,t,x,\mathrm{d}y) \widetilde{P}(t,u,y,B) &=& \int\int P
\bigl(abu,u,(u/s)^{\kappa}(ab)^{\kappa}x,B\bigr) G_{s,t}(\mathrm{d}a)
G_{t,u}(\mathrm{d}b)
\\
&=& \int P\bigl(ru,u,(u/s)^{\kappa}r^{\kappa}x,B\bigr)
G_{s,u}(\mathrm{d}r) = \widetilde {P}(s,u,x,B);
\end{eqnarray*}
in other words, $\widetilde{P}$ satisfies the Chapman--Kolmogorov equations.
\end{pf}

%pr2.2 #&#
\begin{proposition}
If $G_{s,t}$ depends on $s$ and $t$ only through $t/s$ (i.e. $G_{s,t} =
G_{t/s}$), then the scaling property of $P$ carries over to $\widetilde{P}$:
\[
\widetilde{P}\bigl(cs,ct,c^\kappa x,c^\kappa\, \mathrm{d}y\bigr) =
\widetilde{P}(s,t,x,\mathrm{d}y)\quad\quad \forall c>0, s\le t \mbox{ and } x.
\]
\end{proposition}

\begin{pf}
This follows immediately from the definition of $\widetilde{P}$ and the
scaling of $P$.
\end{pf}

Let, for $s\leq t$, $R_{s,t}$ be a random variable having distribution
$G_{s,t}$.
Property \eqref{E:PropG} is equivalent to the property that if, for
$s\le t\le u$, $R_{s,t}$ and $R_{t,u}$ are independent random
variables, then $R_{s,t}R_{t,u} \stackrel{d}{=} R_{s,u}$.
Further, if we let $V_{a,b} = -\ln R_{\mathrm{e}^a,\mathrm{e}^b}$ and write $K_{a,b}$ for
the distributions of $V_{a,b}$ (with $K_{a,a} = \delta_0$),
then Property \eqref{E:PropG} translates to the convolution identity
\[
K_{a,b} \ast K_{b,c} = K_{a,c},\quad\quad a\le b\le c.
\]
As we seek to retain the scaling property of the original process, we
assume that $G_{s,t} = G_{t/s}$ and immediately reduce Property \eqref
{E:PropG} to
\[
K_a \ast K_b = K_{a+b}, \quad\quad a,b\ge0.
\]

The family $(K_a)_{a\geq0}$ defines a subordinator (process with
positive, independent and stationary increments)
and by L\'evy--Khintchine it has Laplace transforms of the form
\[
\int \mathrm{e}^{-\lambda v}K_a(\mathrm{d}v) = \exp\bigl(-a\psi(\lambda)\bigr),
\quad\quad\lambda\ge0.
\]
The function $\psi$, known as the Laplace exponent, takes the form
%
%e2.2 #&#
\begin{equation}
\label{E:LaplaceExp} \psi(\lambda) = \beta\lambda+ \int_0^\infty
\bigl(1-\mathrm{e}^{-\lambda x}\bigr)\nu(\mathrm{d}x),
\end{equation}
with drift $\beta\geq0$ and L\'evy measure $\nu$ satisfying $\nu(\{0\}
)=0$ and $\int_0^\infty(1\wedge x)\nu(\mathrm{d}x)<\infty$.
Conversely, to each $\psi$ (i.e., to each pair $(\beta,\nu)$ as above)
corresponds a convolution semigroup $(K_a)_{a\geq0}$ and in turn a
family $(G_u)_{u\geq1}$ which satisfies Property \eqref{E:PropG}.
For details see, for example, \cite{Ber99},\vspace*{1pt} Section~1.2.
This ensures the existence of $(G_u)_{u\geq1}$ and a process with
transition function $\widetilde{P}$.

%th2.1 #&#
\begin{theorem}\label{T:Transition}
Let $Z$ be a $\kappa$-self-similar Markov process.
To each $\psi$, Laplace exponent of a subordinator, corresponds a
$\kappa$-self-similar Markov process $X$, starting from 0 and having
the marginals of $Z$.
Furthermore, if $Z$ is a martingale and $\psi(\kappa)=\kappa$, then $X$
is also a martingale.

Writing in terms of $R_{t/s}$, $s\le t$, the new transition function
\[
\widetilde{P}(s,t,x,\mathrm{d}y) = \mathbb{E} \bigl[ P\bigl(R_{t/s}t, t,
(t/s)^{\kappa
}R_{t/s}^{\kappa}x, \mathrm{d}y\bigr) \bigr],\quad\quad s\le t
\]
can be seen as a randomisation of $P(s,t,x,\mathrm{d}y)$.
Furthermore, the condition on $(G_{t/s})_{s\leq t}$ for $X$ to be a
martingale can be written as $\mathbb{E}[R_{t/s}^{\kappa}] =
(s/t)^{\kappa}$.
\end{theorem}

\begin{pf}
By\vspace*{1pt} the Kolmogorov extension theorem, there exists a Markov process $X$
with transition function $\widetilde{P}(s,t,x,\mathrm{d}y)$.
As for the martingale property, we first observe that
\[
\psi(\lambda) = -\frac{1}a\ln\int \mathrm{e}^{-\lambda v}K_a(\mathrm{d}v) =
-\frac{1}a\ln\int r^{\lambda}G_{\mathrm{e}^a}(\mathrm{d}r),
\]
so that $\psi(\kappa)=\kappa$ translates to
$u^{\kappa}\int r^{\kappa}G_u(\mathrm{d}r) = 1$ for $u\geq1$.
Then we have
\[
\int y \widetilde{P}(s,t,x,\mathrm{d}y) = \int(t/s)^{\kappa}r^{\kappa}x
G_{t/s}(\mathrm{d}r) = x
\]
using the definition of $\widetilde{P}$ and the martingale property of $Z$.
\end{pf}

The process $X$ can also be obtained using subordination (with a
suitable scaling in the state space and an appropriate time-change).
The idea of subordination was suggested by Bertoin\footnote{Private
communication.} in the context of Brownian marginals. However,
subordination alone (albeit with a suitable scaling in the state space)
is not sufficient. A further logarithmic change of time is needed. This
naturally creates an issue at 0. To deal with this, we follow
Oleszkiewicz \cite{Ole08}.

%pr2.3 #&#
\begin{proposition} \label{T:Timechange}
Let $(\zeta_t)_{t\ge0}$ be a subordinator independent of $Z$.
Let, for $a\in\mathbb{R}$,
\[
X^{(a)}_t = t^{\kappa} \mathrm{e}^{-\kappa\zeta_{a+\ln t}}
Z_{\mathrm{e}^{\zeta_{a+\ln
t}}},\quad\quad t \ge \mathrm{e}^{-a}.
\]
Then the process $(X_t^{(a)})_{t\ge \mathrm{e}^{-a}}$ has the same marginal
distributions as $(Z_t)_{t\ge \mathrm{e}^{-a}}$
and there exists a process $(X_t)_{t\ge0}$ such that for any $a\in
\mathbb{R}$,
$(X_t)_{t\ge \mathrm{e}^{-a}} \stackrel{d}{=} (X^{(a)}_t)_{t\ge \mathrm{e}^{-a}}$.
The process $(X_t)_{t\ge0}$ is a $\kappa$-self-similar Markov process
with transition function
\begin{eqnarray*}
Q(0,t,0,\mathrm{d}y) &=& P(0,t,0,\mathrm{d}y),
\\
Q(s,t,x,\mathrm{d}y) &=& \mathbb{E} \bigl[P\bigl(\mathrm{e}^{-\zeta_{\ln(t/s)}}t,t,(t/s)^{\kappa
}\mathrm{e}^{-{\kappa}\zeta_{\ln(t/s)}}x,\mathrm{d}y
\bigr) \bigr],\quad\quad s\le t.
\end{eqnarray*}
Moreover, $X$ is a martingale provided that $Z$ is a martingale and
$\mathbb{E} [\mathrm{e}^{-{\kappa}\zeta_{\ln(t/s)}} ] = (s/t)^{\kappa}$.
\end{proposition}

\begin{pf}
Since $Z_t(\omega)$ is measurable in $\omega$ and right-continuous in
$t$, it is measurable as a function of $(t,\omega)$.
Hence, for each $t\ge \mathrm{e}^{-a}$, $X_t^{(a)}(\omega)$ is a random variable.

For $c>0$, let $\widehat{Z}_s = \mathrm{e}^{-\kappa\zeta_c}Z_{s\mathrm{e}^{\zeta_c}}$ and
$\widehat{\zeta}_s = \zeta_{c+s} - \zeta_c$, so that $(\widehat
{Z}_s)_{s\ge0} \stackrel{d}{=} (Z_s)_{s\ge0}$ and $(\widehat{\zeta
}_s)_{s\ge0} \stackrel{d}{=} (\zeta_s)_{s\ge0}$.
Then we have, for $t \ge \mathrm{e}^{-b} \ge \mathrm{e}^{-a}$ and with $c=a-b$,
\[
\bigl(X^{(a)}_t \bigr)_{t\ge \mathrm{e}^{-b}} =
\bigl(t^{\kappa} \mathrm{e}^{-{\kappa}\widehat{\zeta}_{b+\ln t}} \widehat {Z}_{\mathrm{e}^{\widehat{\zeta}_{b+\ln t}}}
\bigr)_{t\ge \mathrm{e}^{-b}} \stackrel{d} {=} \bigl(X^{(b)}_t
\bigr)_{t\ge \mathrm{e}^{-b}}.
\]
For $\mathrm{e}^{-a}\le t_1<\cdots<t_n$, let $\mu_{t_1,\ldots,t_n}$ be the law
of $(X^{(a)}_{t_1},\ldots,X^{(a)}_{t_n})$. Then the measures ($\mu
_{t_1,\ldots,t_n})_{n,0<t_1<\cdots<t_n}$ are consistent and by the
Kolmogorov extension theorem, there exists a process $(X_t)_{t>0}$ with
finite-dimensional distributions ($\mu_{t_1,\ldots,t_n})_{n,0<t_1<\cdots<t_n}$.

A similar argument shows that $(X^{(a)}_{ct_1},\ldots,X^{(a)}_{ct_n})
\stackrel{d}{=} c^\kappa(X^{(a)}_{t_1},\ldots,X^{(a)}_{t_n})$, from
which we deduce that $(X_t)_{t>0}$ is $\kappa$-self-similar. As such,
$X$ extends by continuity to $t\geq0$ by letting $X_0=0$.
The equality of marginal distributions of $X^{(a)}$ and $Z$ follows
from the scaling property of $Z$ as
$X^{(a)}_t = t^{\kappa} \mathrm{e}^{-\kappa\zeta_{a+\ln t}} Z_{\mathrm{e}^{\zeta_{a+\ln
t}}} \stackrel{d}{=} t^{\kappa}Z_1 \stackrel{d}{=} Z_t$
for any fixed $t\ge \mathrm{e}^{-a}$.
%$Q(0,t,0,dy) = P(0,t,0,dy)$ for any $t>0$

Using successively Lemmas \ref{L:Lamperti}, \ref{L:GenBochner}, \ref
{L:GenProduct} and \ref{L:GenTimechange} (see the \hyperref[app]{Appendix}),
we see that $X^{(a)}$ is Markovian.
By the scaling property of $P$ and the stationarity of subordinator
$\zeta$, we obtain,
for $\mathrm{e}^{-a}\le s\le t$, the transition function of $X^{(a)}$ as
\[
Q^{(a)}(s,t,x,\mathrm{d}y) = \mathbb{P}\bigl(X^{(a)}_t\in
\mathrm{d}y | X^{(a)}_s=x\bigr) = \mathbb{E} \bigl[P
\bigl(\mathrm{e}^{-\zeta_{\ln(t/s)}}t, t, (t/s)^{\kappa}\mathrm{e}^{-\kappa
\zeta_{\ln(t/s)}}x, \mathrm{d}y\bigr)
\bigr].
\]

Notice that the transition function $Q^{(a)}$ does not depend on $a$
and it is the same as $\widetilde{P}$ defined earlier with $G_{t/s}$
being the distribution of $\exp(-\zeta_{\ln(t/s)})$.
The rest of the assertions of the proposition then follows immediately
from Theorem~\ref{T:Transition}.
(Alternatively, we can carry out the proof independently, without
referring to $\widetilde{P}$, following Oleszkiewicz \cite{Ole08}.)
\end{pf}

The process constructed in Proposition~\ref{T:Timechange} is identical
(in law) to the process obtained in Theorem~\ref{T:Transition} with
$G_{t/s}$ being the distribution of $\exp(-\zeta_{\ln(t/s)})$, or
$R_{t/s} = \exp(-\zeta_{\ln(t/s)})$.

%re2.1 #&#
\begin{remark}
For $X$ to be a martingale, $\beta$ in \eqref{E:LaplaceExp} must be at
most 1, and is 1 if and only if $\zeta_t=t$
($\nu=0$ and $\psi(\lambda)=\lambda$ for any $\lambda$) and $X=Z$.
\end{remark}

%re2.2 #&#
\begin{remark}
In Proposition~\ref{T:Timechange}, we cannot replace $\zeta_{a+\ln t}$
with a two-sided subordinator $\zeta_t = \zeta^1_t \mathbf{1}_{t\ge0} -
\zeta^2_{-t} \mathbf{1}_{t<0}$ for $t\in\mathbb{R}$, where $(\zeta
^1_t)_{t\ge0}$ and $(\zeta^2_t)_{t\ge0}$ are independent subordinators.
This is because by doing that, we will not have independent increments.
In particular, since $\zeta_0=0$, then for $t<0$, the increment $\zeta
_0 - \zeta_t = -\zeta_t \in\mathcal{G}_t$, where $\mathcal{G}_t$
denotes the filtration generated by $\zeta$.
\end{remark}

%s3 #&#
\section{Properties}

In this section, we obtain the infinitesimal generators of the process
$X$ and display some of their path properties.
We will work within the martingale framework, that is, unless otherwise
stated, we will assume that our initial process $Z$ is a martingale and
we will use a subordinator $\zeta$
with drift $\beta$, L\'evy measure $\nu$ and Lapace exponent satisfying
$\psi(\kappa)=\kappa$.

%th3.1 #&#
\begin{theorem} \label{T:Generator}
Suppose that $Z$ has infinitesimal generator $L_t$.
Then the infinitesimal generator of the process $X$ is given by
\begin{eqnarray*}
A_0f(x) &=& L_0f(x),
\\
A_tf(x) &=& \beta L_t f(x) + (1-\beta) \frac{\kappa}{t}
x f'(x)
\\
&&{} + \frac{1}{t} \int_{(0,\infty)} \int_{-\infty}^\infty
\bigl(f(y)-f(x)\bigr) P\bigl(t\mathrm{e}^{-u}, t, x\mathrm{e}^{-u\kappa},\mathrm{d}y\bigr)
\nu(\mathrm{d}u),\quad\quad t>0,
\end{eqnarray*}
for $f$ differentiable and in the domain of $L$.
\end{theorem}

\begin{pf}
First, from Lemma~\ref{L:Lamperti}, $\widehat{Z}_t := \mathrm{e}^{-t\kappa}
Z_{\mathrm{e}^t}$ is time-homogeneous with generator
$\widehat{L}f(x) = L_1 f(x) - {\kappa}xf'(x)$
and transition semigroup
$\widehat{P}_tf(x) = \int f(y) P(\mathrm{e}^{-t},1,x\mathrm{e}^{-t\kappa},\mathrm{d}y)$.
Next, applying Lemma~\ref{L:GenBochner}, the generator of the process
$\bar{Z}_t := \widehat{Z}_{\zeta_t} = \mathrm{e}^{-{\kappa}\zeta_t}Z_{\mathrm{e}^{\zeta
_t}}$ is
\[
\bar{L}f(x) = \beta L_1 f(x) - {\kappa}\beta xf'(x) +
\int_{(0,\infty)} \int\bigl(f(y)-f(x)\bigr) P\bigl(\mathrm{e}^{-u},1,
x\mathrm{e}^{-u\kappa},\mathrm{d}y
\bigr) \nu(\mathrm{d}u).
\]
Then, let $\widetilde{Z}_t = \mathrm{e}^{{\kappa}(t-a)}\bar{Z}_t = \mathrm{e}^{{\kappa
}(t-a)}\mathrm{e}^{-{\kappa}\zeta_t}Z_{\mathrm{e}^{\zeta_t}}$
and using Lemma~\ref{L:GenProduct}, the generator of $\widetilde{Z}$ is
\begin{eqnarray*}
\widetilde{L}_tf(x) &=& \beta\pi_{\mathrm{e}^{-{\kappa}(t-a)}} L_1
\pi_{\mathrm{e}^{{\kappa}(t-a)}} f(x) + (1-\beta) {\kappa}xf'(x)
\\
&&{} + \int_{(0,\infty)} \int \bigl(f\bigl(\mathrm{e}^{{\kappa}(t-a)}y\bigr) -
f(x) \bigr) P\bigl(\mathrm{e}^{-u},1, \mathrm{e}^
{-{\kappa}(t-a)}x\mathrm{e}^{-u\kappa},\mathrm{d}y\bigr)
\nu(\mathrm{d}u),
\end{eqnarray*}
since $\pi_{\mathrm{e}^{-{\kappa}(t-a)}} \Lambda\pi_{\mathrm{e}^{{\kappa}(t-a)}} f(x) =
\Lambda f(x)$ for $\Lambda f(x) = xf'(x)$.
Finally, we time-change the process $\widetilde{Z}_t$ with $a+\ln t$ to
get $X^{(a)}_t$.
Thus, by Lemma~\ref{L:GenTimechange}, the generator of $X^{(a)}_t$ is
\begin{eqnarray*}
A^{(a)}_tf(x) &=& \beta L_t f(x) + (1-\beta)
\frac{\kappa}{t} xf'(x)
\\
&&{} + \frac{1}{t} \int_{(0,\infty)} \int \bigl(f(y) - f(x)
\bigr) P\bigl(t\mathrm{e}^{-u}, t, x\mathrm{e}^{-u\kappa}, \mathrm{d}y\bigr) \nu(\mathrm{d}u)
\end{eqnarray*}
due to a change of variable, the scaling property of $P$ and identity
\eqref{E:SSGen1}.
The generator of $X_t$ is established by noting that $A^{(a)}_t$ does
not depend on $a$.

Since $Z$ is self-similar, $\pi_cf$ is in the domain of $L$ for all
$c>0$ whenever $f$ is in the domain of $L$.
Therefore, $f$ is in the domain of $A$, if $f$ is also differentiable.
\end{pf}

Note that when $\beta=1$ and $\nu\equiv0$, we recover the process $Z$
and $A_t=L_t$.

For a measurable function $f$, if there exists a measurable function
$g$ such that
for each $t$, $\int_0^t |g(X_s)| \,\mathrm{d}s <\infty$ almost surely
and $f(X_t) - f(X_0) - \int_0^t g(X_s) \,\mathrm{d}s$ is a martingale,
then $f$ is said to belong to the domain of the extended infinitesimal
generator of $X$ and
the extended infinitesimal generator $A_sf(X_s) = g(X_s)$.
If $f(x)=x^2$ belongs to the domain of the extended infinitesimal
generator of $X$,
then $X$ has predictable quadratic variation
%
%e3.1 #&#
\begin{equation}
\label{E:QuadVar} \anglel X,X\angler _t = \int_0^t
A_sf(X_{s})\,\mathrm{d}s.
\end{equation}
See examples in Section~\ref{S:Examples} for the computation $\anglel X,X\angler $ in some specific cases.

%pr3.1 #&#
\begin{proposition}
Suppose that $Z$ is continuous in probability.
Then the process $X$ is also continuous in probability, that is, for
every $t$,
\[
\forall c>0,\quad\quad \lim_{s\to t} \mathbb{P} \bigl(\llvert
X_t - X_s \rrvert >c \bigr) = 0.
\]
\end{proposition}

\begin{pf}
We have, for $s,t \ge \mathrm{e}^{-a}$,
%
%e3.2 #&#
\begin{eqnarray}\label{E:ContProb1}
\mathbb{P} \bigl(\llvert X_t - X_s \rrvert >c \bigr) &
\le&\mathbb{P} \biggl( \bigl| \bigl(t^{\kappa} \mathrm{e}^{-{\kappa}\zeta_{a+\ln t}} - s^{\kappa}
\mathrm{e}^{-{\kappa}\zeta_{a+\ln s}}\bigr) Z_{\mathrm{e}^{\zeta_{a+\ln t}}} \bigr| > \frac{c}{2} \biggr)
\nonumber
\\[-8pt]\\[-8pt]
&&{} + \mathbb{P} \biggl( \bigl| s^{\kappa} \mathrm{e}^{-{\kappa}\zeta_{a+\ln s}} (Z_{\mathrm{e}^{\zeta_{a+\ln t}}} -
Z_{\mathrm{e}^{\zeta_{a+\ln s}}}) \bigr| > \frac{c}{2} \biggr).\nonumber
\end{eqnarray}
However, the first term
\begin{eqnarray*}
&&\mathbb{P} \biggl( \bigl| \bigl(t^{\kappa} \mathrm{e}^{-{\kappa}\zeta_{a+\ln t}} - s^{\kappa}
\mathrm{e}^{-{\kappa}\zeta_{a+\ln s}}\bigr) Z_{\mathrm{e}^{\zeta_{a+\ln t}}} \bigr| > \frac{c}{2} \biggr)
\\
&&\quad\le\mathbb{P} \biggl( \bigl| \bigl(t^{\kappa} - s^{\kappa}\bigr)
\mathrm{e}^{-{\kappa}\zeta
_{a+\ln t}} Z_{\mathrm{e}^{\zeta_{a+\ln t}}} \bigr| > \frac{c}{4} \biggr)
\\
&&\quad\quad{}+ \mathbb{P}
\biggl( \bigl| s^{\kappa} \bigl(\mathrm{e}^{-{\kappa}\zeta_{a+\ln t}} - \mathrm{e}^{-{\kappa}\zeta_{a+\ln s}}\bigr)
Z_{\mathrm{e}^{\zeta_{a+\ln t}}} \bigr| > \frac
{c}{4} \biggr),
\end{eqnarray*}
which converges to 0 as $s\to t$, since $\zeta$ is continuous in
probability as a subordinator.

To deal with the last term in \eqref{E:ContProb1} we first observe that a
process that is continuous in probability does not jump at fixed points
so that $\mathbb{P}(\Delta \mathrm{e}^{\zeta_{a+\ln t}} = 0) = 1$.
Further,\vadjust{\goodbreak} since $Z$ is also continuous in probability,
$\lim_{s\to t} \mathbb{P} (  | Z_{y_t} - Z_{y_s}  | >
\varepsilon ) = 0$ as soon as $\lim_{s\to t} y_s = y_t$.
Therefore, for $s\leq t+1$,
\begin{eqnarray*}
&&\lim_{s\to t}  \mathbb{P} \biggl( \bigl| s^{\kappa}
\mathrm{e}^{-{\kappa}\zeta
_{a+\ln s}} (Z_{\mathrm{e}^{\zeta_{a+\ln t}}} - Z_{\mathrm{e}^{\zeta_{a+\ln s}}}) \bigr| > \frac{c}{2}
\biggr)
\\
&&\quad\le\lim_{s\to t} \mathbb{P} \biggl(  | Z_{\mathrm{e}^{\zeta_{a+\ln t}}} -
Z_{\mathrm{e}^{\zeta_{a+\ln s}}}  | > \frac{c}{2} (t+1)^{-\kappa} \biggr)
\\
&&\quad= \mathbb{E} \biggl[ \lim_{s\to t} \mathbb{P} \biggl(
 |Z_{\mathrm{e}^{\zeta
_{a+\ln t}}} - Z_{\mathrm{e}^{\zeta_{a+\ln s}}}  | > \frac{c}{2}(t+1)^{-\kappa
}
\big| \hat{\mathcal{G}}_t \biggr) \mathbf{1}_{ \{ \Delta \mathrm{e}^{\zeta
_{a+\ln t}} = 0  \}} \biggr] = 0,
\end{eqnarray*}
where $\hat{\mathcal{G}}_t = \sigma(\zeta_s, s\le a+\ln t)$.
\end{pf}

%pr3.2 #&#
\begin{proposition}
If $Z$ is continuous in probability with finite second moments and
$\zeta$ has no drift $(\beta=0)$, then $X$ is a purely discontinuous martingale.
\end{proposition}

\begin{pf}
Let $U_t = \mathrm{e}^{\zeta_{a+\ln t}}$ so that $X_t = t^{\kappa} U_t^{-\kappa
}Z_{U_t}$ for $t\ge \mathrm{e}^{-a}$.
First, we observe that with probability one, $Z$ does not jump at
$U_{t-}$ if $U$ jumps at $t$.
Indeed, if $\Gamma$ is a countable set of points in $[0,\infty)$, then,
as $Z$ is continuous in probability, $\mathbb{P}(\exists t \in\Gamma
\mbox{ s.t. } \Delta Z_t \neq0) = 0$.
Let $\Lambda= \{t>\mathrm{e}^{-a}\dvt \Delta U_t >0\}$ and $\Gamma_U = U^-(\Lambda
)$ where $U^-_t=U_{t-}$, then
\[
\mathbb{P} \bigl(\exists t>\mathrm{e}^{-a} \mbox{ s.t. } \Delta
U_t>0, \Delta Z_{(U_{t-})}\neq0\bigr) = \mathbb{E}\bigl[
\mathbb{P} (\exists s\in\Gamma_U \mbox{ s.t. } \Delta
Z_s\neq0 | \mathcal{G}_\infty)\bigr] = 0,
\]
where $\mathcal{G}$ denotes the filtration of $U$.
Taking $a$ to infinity, we obtain the desired result.

To show that $X$ is purely discontinuous, that is, $\anglel X^c,X^c\angler _t = 0$, we compute the sum of the square of jumps of $X$. In general,
we have
\begin{eqnarray*}
\mathbb{E} \biggl[ \sum_{\mathrm{e}^{-a}<s\le t} (\Delta
X_s)^2 \biggr] &=& \mathbb {E} \biggl[ \sum
_{\mathrm{e}^{-a}<s\le t} (\Delta X_s)^2
\mathbf{1}_{\Delta U_s
>0, \Delta Z_{(U_{s-})}\neq0} \biggr]
\\
&&{}+ \mathbb{E} \biggl[ \sum_{\mathrm{e}^{-a}<s\le t} (\Delta
X_s)^2 \mathbf {1}_{\Delta U_s >0, \Delta Z_{(U_{s-})}=0} \biggr]
\\
&&{}+ \mathbb{E}
\biggl[ \sum_{\mathrm{e}^{-a}<s\le t} (\Delta X_s)^2
\mathbf{1}_{\Delta U_s =0} \biggr].
\end{eqnarray*}
As $Z$ is continuous in probability, the first term is zero due to the
observation at the start of the proof.
We write $l(a, t) = \mathbb{E} [ \sum_{\mathrm{e}^{-a}<s\le t} (\Delta X_s)^2
\mathbf{1}_{\Delta U_s =0}  ]$ for the third term.

For the second term, we have, on the set $\{\Delta U_s>0, \Delta
Z_{(U_{s-})}=0\}$,
\begin{eqnarray*}
(\Delta X_s)^2 &=& s^{2\kappa} \bigl(
U_s^{-\kappa}(Z_{U_s}-Z_{(U_{s-})}) +
Z_{(U_{s-})}\bigl(U_s^{-\kappa}-U_{s-}^{-\kappa}
\bigr) \bigr)^2
\\
&=& s^{2\kappa} \bigl( U_s^{-2\kappa}(Z_{U_s}-Z_{(U_{s-})})^2
+ Z_{(U_{s-})}^2\bigl(U_s^{-\kappa}-U_{s-}^{-\kappa}
\bigr)^2
\\
&&\hphantom{s^{2\kappa} \bigl(}{} + 2U_s^{-\kappa}\bigl(U_s^{-\kappa}-U_{s-}^{-\kappa
}
\bigr)Z_{(U_{s-})}(Z_{U_s}-Z_{(U_{s-})}) \bigr).
\end{eqnarray*}
Let $\theta_t = \mathbb{E}[Z_t^2]$. Since $Z$ is $\kappa$-self-similar,
$\theta_t = t^{2\kappa}\theta_1$. As $Z$ is a martingale, $\mathbb
{E}[Z_s(Z_t-Z_s)]=0$ and $\mathbb{E}[(Z_t-Z_s)^2]=\mathbb
{E}[Z_t^2]-\mathbb{E}[Z_s^2]$. Thus, we obtain
\begin{eqnarray*}
&&\mathbb{E}\bigl[(\Delta X_s)^2 | \mathcal{G}_\infty
\cap\{\Delta U_s>0, \Delta Z_{(U_{s-})}=0\}\bigr]
\\
& &\quad= s^{2\kappa} \bigl( U_s^{-2\kappa}(
\theta_{U_s}-\theta_{(U_{s-})}) + \theta_{(U_{s-})}
\bigl(U_s^{-\kappa}-U_{s-}^{-\kappa}
\bigr)^2 \bigr)
\\
&&\quad = 2 s^{2\kappa} \theta_1 \bigl( 1 - U_s^{-\kappa}U_{s-}^{\kappa}
\bigr).
\end{eqnarray*}
Since $\{\Delta Z_{(U_{s-})}=0\}$ has probability one on the set $\{
\Delta U_s >0\}$, it follows that
\begin{eqnarray*}
\mathbb{E} \biggl[ \sum_{\mathrm{e}^{-a}<s\le t} (\Delta
X_s)^2 \mathbf{1}_{\Delta
U_s >0, \Delta Z_{(U_{s-})}=0} \biggr] &=& \mathbb{E}
\biggl[ \sum_{\mathrm{e}^{-a}<s\le t} 2 s^{2\kappa}
\theta_1 \bigl( 1 - U_s^{-\kappa}U_{s-}^{\kappa}
\bigr) \mathbf{1}_{\Delta U_s>0} \biggr]
\\
&=& \mathbb{E} \biggl[ \sum_{0<r\le a+\ln t} 2 \mathrm{e}^{2\kappa(r-a)}
\theta_1 \bigl( 1 - \mathrm{e}^{-{\kappa}\Delta\zeta_r} \bigr) \mathbf{1}_{\Delta\zeta
_r>0}
\biggr].
\end{eqnarray*}
Writing in terms of $\nu$ and $\psi$,
and using \eqref{E:LaplaceExp} with $\psi(\kappa) = \kappa$,
\begin{eqnarray*}
\mathbb{E} \biggl[ \sum_{\mathrm{e}^{-a}<s\le t} (\Delta
X_s)^2 \mathbf{1}_{\Delta
U_s >0,\Delta Z_{(U_{s-})}=0} \biggr] &=&
\theta_1 \int_0^{a+\ln t} 2
\mathrm{e}^{2\kappa(r-a)} \,\mathrm{d}r \int_{(0,\infty)}
\bigl(1-\mathrm{e}^{-z\kappa}\bigr)
\nu(\mathrm{d}z)
\\
&=& \theta_1 \bigl(t^{2\kappa}-\mathrm{e}^{-2a\kappa}\bigr) (1-\beta).
\end{eqnarray*}

Adding those three terms and taking limit as $a\to\infty$, we have
\[
\mathbb{E} \biggl[ \sum_{0<s\le t}  (\Delta
X_s)^2 \biggr] = \theta_1 t^{2\kappa}
(1-\beta) + \lim_{a\to\infty} l(a,t).
\]
Since $X$ is square integrable on any finite interval if $Z$ has finite
second moments,
it has quadratic variation with expectation
$\mathbb{E} [ [X,X]_t  ] = \mathbb{E} [X_t^2 ] = \mathbb
{E} [Z_t^2 ] = t^{2\kappa} \theta_1$.
Therefore,
\begin{eqnarray*}
\mathbb{E} \bigl[ \bigl\anglel X^c,X^c\bigr\angler
_t \bigr] &=& \mathbb{E} \bigl[ [X,X]_t \bigr] -
\mathbb{E} \biggl[  \sum_{0<s\le t, \Delta X_s\neq0}  (\Delta
X_s)^2 \biggr]
\\
&=& \theta_1 t^{2\kappa} \beta- \lim_{a\to\infty}
l(a,t).
\end{eqnarray*}
Note that both $\mathbb{E} [ \anglel X^c,X^c\angler _t  ]$ and
$l(a,t)$ are non-negative.
Thus, when $\beta=0$, we must have $\lim_{a\to\infty} l(a,t) = 0$,
which gives $\mathbb{E} [ \anglel X^c,X^c\angler _t  ] = 0$ and
hence $\anglel X^c,X^c\angler _t = 0$.
\end{pf}

%s4 #&#
\section{Examples} \label{S:Examples}

Given any self-similar Markov martingale $Z$ with transition function
$P$ and infinitesimal generator $L$, we can mimic $Z$ as per Section~\ref{S:Scheme}.
We construct a new Markov martingale $X$ that has the same marginal
distributions as $Z$
and possesses the same self-similarity\vadjust{\goodbreak} $Z$ enjoys from each~$\zeta$.
We assume that the subordinator $\zeta$ has drift $\beta$ and L\'evy
measure $\nu$ with Laplace exponent $\psi$ satisfying $\psi(\kappa
)=\kappa$, or $\mathbb{E} [\mathrm{e}^{-\kappa\zeta_{\ln(t/s)}} ] =
(s/t)^\kappa$.

Examples of subordinators include Poisson process, compound Poisson
process with positive jumps,
gamma process and stable subordinators.
For example, we can take $\zeta$ as a Poisson process with rate $\kappa
/(1-\exp(-\kappa))$ to satisfy $\psi(\kappa) = \kappa$.
In the following, we provide some examples of mimicking with the
infinitesimal generators and the predictable quadratic variations
computed explicitly to have a better understanding of the processes.

We finish this section with a discussion on modifying our construction
to mimic some Brownian related martingales and its limitation.

%s4.1 #&#
\subsection{Gaussian continuous martingales} \label{S:GausContMart}

For any $k\ge0$, define the process $(Z_t)_{t\ge0}$ by $Z_t = \int_0^t
s^k \,\mathrm{d}B_s$,
where $(B_t)_{t\ge0}$ is a Brownian motion.
Note that
\[
(Z_t)_{t\ge0} \stackrel{d} {=} (B_{(1/(2k+1))t^{2k+1}}
)_{t\ge0}.
\]
This is a Gaussian process with zero mean and covariance function
$\operatorname{Cov}(Z_t,Z_{t+u}) = \frac{1}{2k+1}t^{2k+1}$.
It is a Markov process and also a martingale.
Moreover, it is $(k+\frac{1}{2})$-self-similar ($\kappa=k+\frac{1}2$)
since for all $c>0$,
\[
(Z_{ct})_{t\ge0} \stackrel{d} {=} (B_{(1/(2k+1))c^{2k+1}t^{2k+1}}
)_{t\ge0} \stackrel{d} {=} \bigl(c^{k+ 1/2}B_{(1/(2k+1))t^{2k+1}}
\bigr)_{t\ge0} \stackrel{d} {=} \bigl(c^{k+1/2}Z_t
\bigr)_{t\ge0}.
\]

A key aspect of the construction in \cite{HamKle07} is the following
representation of the mimic $X_t$ when $k=0$:
\[
X_t = \sqrt{t/s} ( \sqrt{R_{s,t}} X_s +
\sqrt{s} \sqrt{1-R_{s,t}} \xi_{s,t} ),\quad\quad t\ge s,
\]
where $R_{s,t}$ has distribution $G_{t/s}$, $\xi_{s,t}$ is standard
normal and, $R_{s,t}$, $\xi_{s,t}$ and $X_s$ are independent.
This representation extends to the case of other Gaussian continuous
martingales.
In fact, it also extends to the case of stable processes -- see
Proposition~\ref{P:StableRep}.

%pr4.1 #&#
\begin{proposition} \label{P:GausContMartRep}
With $\kappa= k+\frac{1}{2}$, the mimic $(X_t)_{t\ge0}$ has the representation
\[
X_t = (t/s )^{\kappa} \bigl( R_{s,t}^\kappa
X_s + s^{\kappa} (1-R_{s,t})^{\kappa}
\xi_{s,t} \bigr),\quad\quad t\ge s>0,
\]
where $R_{s,t}\stackrel{d}{=}\mathrm{e}^{-\zeta_{\ln(t/s)}}$, $\xi
_{s,t}\stackrel{d}{=}Z_1$ and, $R_{s,t}$, $\xi_{s,t}$ and $X_s$ are independent.
\end{proposition}

\begin{pf}
Since $Z$ and $\zeta$ have independent and stationary increments, for
$\mathrm{e}^{-a} \le s < t$,
\begin{eqnarray*}
X_t^{(a)} &\stackrel{d} {=}& t^{\kappa}
\mathrm{e}^{-{\kappa}\zeta_{a+\ln t}} Z_{\mathrm{e}^{\zeta
_{a+\ln s}}} + t^{\kappa} \mathrm{e}^{-{\kappa}\zeta_{a+\ln t}}
\bigl(\mathrm{e}^{\zeta_{a+\ln
t}} - \mathrm{e}^{\zeta_{a+\ln s}}\bigr)^{\kappa} \xi_{s,t}
\\
&\stackrel{d} {=}& (t/s)^{\kappa} \mathrm{e}^{-{\kappa}\widehat{\zeta}_{\ln
(t/s)}} X_s^{(a)}
+ t^{\kappa} \bigl(1- \mathrm{e}^{-\widehat{\zeta}_{\ln
(t/s)}}\bigr)^{\kappa}
\xi_{s,t},
\end{eqnarray*}
where $\widehat{\zeta}$ is an independent copy of $\zeta$, and $\xi_{s,t}$ is a random variable
distributed as $Z_1$.
Note that this representation holds also for $t=s$.
\end{pf}

Knowing that $Z$ has generator $L_tf(x) = \frac{1}{2}t^{2k}f''(x)$ for
$t\ge0$, we can compute the generator of $X$ following Theorem~\ref
{T:Generator} and obtain, for $t>0$,
\begin{eqnarray*}
A_tf(x) &=& \frac{1}{2}\beta t^{2k}f''(x)
+ (1-\beta)\frac
{2k+1}{2t}xf'(x)
\\
&&{} + \frac{1}{t}\int_{(0,\infty)} \int_{-\infty}^\infty
\bigl(f(y)-f(x)\bigr) P\bigl(t\mathrm{e}^{-u},t,x\mathrm{e}^{-(k+1/2)u},\mathrm{d}y\bigr)
\nu(\mathrm{d}u).
\end{eqnarray*}
Taking $f(x)=x^2$, then Equation \eqref{E:QuadVar}, along with routine
calculations and Equation \eqref{E:LaplaceExp}, gives the following result.

%pr4.2 #&#
\begin{proposition}
The predictable quadratic variation of $X$ is
\[
\anglel X,X\angler _t = \frac{1}{(2k+1)^2}t^{2k+1}\psi(2k+1)
+ \bigl(2k+1-\psi (2k+1)\bigr) \int_0^t
\frac{1}{s}X_s^2 \,\mathrm{d}s.
\]
\end{proposition}

Since $\anglel Z,Z\angler _t = \int_0^t s^{2k}\, \mathrm{d}s$, we can also write
\[
\anglel X,X\angler _t = \frac{1}{2k+1}\psi(2k+1)\anglel Z,Z
\angler _t + \bigl(2k+1-\psi(2k+1)\bigr) \int_0^t
\frac{1}{s} X_s^2 \,\mathrm{d}s.
\]

%re4.1 #&#
\begin{remark}
When $k=0$, $Z$ is a Brownian motion and our results agree with \cite{HamKle07}.
\end{remark}

%s4.2 #&#
\subsection{Martingale of squared Bessel process}

A process $(S_t)_{t\ge0}$ is a squared Bessel process of dimension
$\delta$, for some $\delta\ge0$, if it satisfies
$\mathrm{d}S_t = 2\sqrt{S_t} \,\mathrm{d}B_t + \delta \,\mathrm{d}t$,
where $B$ denotes a Brownian motion.
The squared Bessel process $S$ started at 0 is a continuous Markov
process satisfying the self-similarity with $\kappa=1$.

Let $Z_t = S_t -\delta t$. Then $(Z_t)_{t\ge0}$ is a 1-self-similar
Markov process and satisfies the SDE
\[
\mathrm{d}Z_t = 2\sqrt{Z_t +\delta t} \,\mathrm{d}B_t.
\]
Note that $S$ is stochastically dominated by the square of the norm of
an $n$-dimensional Brownian motion, where $n\geq\delta$, thus $\mathbb
{E}[S_t]\leq nt$ and $\mathbb{E}[\anglel Z,Z\angler _t]\leq2nt^2$. It
follows that $Z$ is a true martingale (see, e.g., \cite{Kle12}, Theorem~7.35).

The infinitesimal generator of $Z$ is
$L_tf(x) = 2(x+\delta t)f''(x)$, $t\ge0$,
thus, that of $X$ is
\begin{eqnarray*}
A_0f(x) &=& L_0f(x) = 2xf''(x),
\\
A_tf(x) &=& 2\beta(x+\delta t)f''(x) +
\frac{1}{t}(1-\beta)xf'(x)
\\
&&{} + \frac{1}{t} \int_{(0,\infty)} \int_{-\infty}^\infty
\bigl(f(y)-f(x)\bigr) P\bigl(t\mathrm{e}^{-u},t,x\mathrm{e}^{-u},\mathrm{d}y\bigr)
\nu(\mathrm{d}u),\quad\quad t>0.
\end{eqnarray*}

%pr4.3 #&#
\begin{proposition}
The predictable quadratic variation of $X$ is
\[
\anglel X,X\angler _t = \delta t^2\psi(2) + 4\bigl(
\psi(2)-1\bigr) \int_0^t X_{s} \,\mathrm{d}s +
\bigl(2-\psi(2)\bigr) \int_0^t
\frac{1}{s}X_{s}^2 \,\mathrm{d}s,\quad\quad t\ge0,
\]
where $\psi$ is the Laplace exponent of $\zeta$.
\end{proposition}

\begin{pf}
From $\mathrm{d}Z_t = 2\sqrt{Z_t +\delta t} \,\mathrm{d}B_t$, we have for $u>0$,
\[
Z_t^2 = Z_{t\mathrm{e}^{-u}}^2 + 4 \int
_{t\mathrm{e}^{-u}}^t Z_s\sqrt{Z_s+
\delta s} \,\mathrm{d}B_s + 4 \int_{t\mathrm{e}^{-u}}^t
Z_s \,\mathrm{d}s + 2 \delta t^2\bigl(1-\mathrm{e}^{-2u}\bigr)
\]
and taking conditional expectation ($\int_0^t Z_s\sqrt{Z_s+\delta s}
\,\mathrm{d}B_s$ is a true martingale -- see above),
\[
\mathbb{E}\bigl[Z_t^2 | Z_{t\mathrm{e}^{-u}}\bigr] =
Z_{t\mathrm{e}^{-u}}^2 + 4t\bigl(1-\mathrm{e}^{-u}\bigr)Z_{t\mathrm{e}^{-u}}
+ 2\delta t^2\bigl(1-\mathrm{e}^{-2u}\bigr).
\]
Thus, using Equation \eqref{E:LaplaceExp} we obtain, with $f(x)=x^2$,
\[
A_tf(x)  = 2\delta t \psi(2) + 4x\bigl(\psi(2)-\psi(1)\bigr) +
\frac
{1}{t}x^2\bigl(2-\psi(2)\bigr).
\]
However $\psi(1)=1$.
The result then follows from Equation \eqref{E:QuadVar}.
Note that if $\beta=1$, $\psi(\lambda)=\lambda$ for any $\lambda$ and
we recover $\anglel X,X\angler _t = 2\delta t^2 + 4\int_0^t X_{s} \,\mathrm{d}s$.
\end{pf}

%s4.3 #&#
\subsection{Stable processes with \texorpdfstring{$1<\alpha<2$}{1<alpha<2}}

Suppose $(Z_t)_{t\ge0}$ is an $\alpha$-stable process with $1<\alpha<2$.
Then $Z$ is a Markov process and
\[
(Z_{c^\alpha t})_{t\ge0} \stackrel{d} {=} (cZ_t)_{t\ge0},\quad\quad
\forall c>0,
\]
that is, $Z$ is $\kappa$-self-similar with $\kappa=\frac{1}{\alpha}$.
It is a L\'evy process with L\'evy triplet $(0,\nu_Z,\gamma)$, where
\[
\nu_Z(\mathrm{d}z) = (A\mathbf{1}_{z>0} + B\mathbf{1}_{z<0})|z|^{-(\alpha+1)}\,\mathrm{d}z
\]
for some positive constants $A$ and $B$.
Assume that $Z$ is a martingale, in which case the L\'evy triplet must satisfy
\[
\gamma+ \int_{|z|\ge1}z\nu_Z(\mathrm{d}z) = 0.
\]

%pr4.4 #&#
\begin{proposition}\label{P:StableRep}
The mimic $(X_t)_{t\ge0}$ has the representation
\[
X_t = (t/s )^{\kappa} \bigl( R_{s,t}^\kappa
X_s + s^{\kappa} (1-R_{s,t})^{\kappa}
\xi_{s,t} \bigr), \quad\quad t\ge s,
\]
where $R_{s,t}\stackrel{d}{=}\mathrm{e}^{-\zeta_{\ln(t/s)}}$, $\xi
_{s,t}\stackrel{d}{=}Z_1$ and, $R_{s,t}$, $\xi_{s,t}$ and $X_s$ are independent.
\end{proposition}

\begin{pf}
See Proposition~\ref{P:GausContMartRep}.
\end{pf}

The stable process $Z$ has infinitesimal generator
\[
Lf(x) = \gamma f'(x) + \int_{\mathbb{R}\setminus\{0\}} \bigl(f(x+y)
- f(x) - yf'(x)\mathbf{1}_D(y) \bigr)
\nu_Z(\mathrm{d}y),
\]
where $D=\{x\dvt|x|\le1\}$.
To distinguish the characteristics of $\zeta$ from that of $Z$, we add
the subscript $\zeta$ to the drift and L\'evy measure of the
subordinator $\zeta$.
Then Theorem~\ref{T:Generator} gives the infinitesimal generator of $X$
for $t>0$ as
\begin{eqnarray*}
A_tf(x) &=& \beta_\zeta\gamma f'(x) + (1-
\beta_\zeta) \frac{x\kappa}{t} f'(x)
\\
&&{} + \beta_\zeta\int_{\mathbb{R}\setminus\{0\}} \bigl(f(x+y) - f(x) -
yf'(x)\mathbf{1}_D(y) \bigr) \nu_Z(\mathrm{d}y)
\\
&&{} + \frac{1}{t} \int_{(0,\infty)} \int_{-\infty}^\infty
\bigl(f(w)-f(x)\bigr) P\bigl(t\mathrm{e}^{-u}, t, x\mathrm{e}^{-u\kappa},\mathrm{d}w\bigr)
\nu_\zeta(\mathrm{d}u).
\end{eqnarray*}

%s4.4 #&#
\subsection{Martingale with marginals \texorpdfstring{$t^\kappa V$}{tkappaV} with $V$ symmetric}

Let $V$ be an integrable, symmetric random variable (i.e. $V \stackrel
{d}{=} -V$) and
$(B_t)_{t\ge0}$ be a Brownian motion independent of $V$. Following and
extending \cite{HirProRoyYor11}, page 283, for any $\kappa$ let $T_t =
\inf\{u\ge0 \dvt |B_u| = t^\kappa\}$ and
$Z_t = B_{T_t}V$. Then $(Z_t)_{t\ge0}$ is a Markov martingale such that
for each $t\ge0$, $Z_t \stackrel{d}{=} t^\kappa V$.
Moreover, $Z$ is $\kappa$-self-similar. Indeed, using the Brownian
motion $B_t^{(c)} := cB_{c^{-2}t}$, we have
\[
T_t^{(c)}  = \inf\bigl\{u\ge0 \dvt \bigl|B_u^{(c)}\bigr|
= t^\kappa\bigr\} = c^2\inf\bigl\{s\ge0 \dvt |B_s|
= c^{-1}t^\kappa\bigr\} = c^2 T_{c^{-1/\kappa}t}.
\]
It follows that $B^{(c)}_{T^{(c)}_t} = cB_{c^{-2}c^2T_{c^{-1/\kappa}t}}
= cB_{T_{c^{-1/\kappa}t}}$
and $(B_{T_{ct}})_{t\ge0} \stackrel{d}{=} (c^\kappa B_{T_t})_{t\ge0}$. Hence,
\[
(Z_{ct})_{t\ge0} = (B_{T_{ct}}V)_{t\ge0}
\stackrel{d} {=} \bigl(c^\kappa B_{T_t}V\bigr)_{t\ge0} =
\bigl(c^\kappa Z_t\bigr)_{t\ge0}.
\]

Since $Z$ is a Markov process with transition semigroup
\begin{eqnarray*}
P_{0,t}f(x) &=& \int f\bigl(t^\kappa v\bigr) \,\mathrm{d}F(v),\quad\quad t>0,
\\
P_{s,t}f(x) &=& \frac{1}{2} \bigl(1+(s/t)^\kappa \bigr)
f \bigl((t/s)^\kappa x \bigr) + \frac{1}{2} \bigl(1-(s/t)^\kappa
\bigr) f \bigl(-(t/s)^\kappa x \bigr),\quad\quad 0<s\le t,
\end{eqnarray*}
where $F$ is the cumulative distribution function of $V$, it has
infinitesimal generator
\begin{eqnarray*}
L_0 f(0) &=& 0,
\\
L_t f(x) &=& \frac{\kappa}{t} xf'(x) -
\frac{\kappa}{2t}f(x) + \frac
{\kappa}{2t} f(-x),\quad\quad t>0.
\end{eqnarray*}
Using Theorem~\ref{T:Generator}, we obtain the infinitesimal generator
of the mimic $X$, for $t>0$,
\begin{eqnarray*}
A_t f(x) &=& \frac{\kappa}{t} \biggl( xf'(x) +
\frac{1}{2}\beta f(-x) - \frac{1}{2}\beta f(x) \biggr)
\\
&&{} + \frac{1}{t} \int_{(0,\infty)} \int_{-\infty}^\infty
\bigl(f(y)-f(x)\bigr) P\bigl(t\mathrm{e}^{-u},t,x\mathrm{e}^{-\kappa u},\mathrm{d}y\bigr)
\nu(\mathrm{d}u).
\end{eqnarray*}

%pr4.5 #&#
\begin{proposition}
The predictable quadratic variation of $X$ is
\[
\anglel X,X \angler _t = 2\kappa\int_0^t
\frac{1}{s}X_{s}^2\, \mathrm{d}s,\quad\quad t\ge0.
\]
\end{proposition}

\begin{pf}
Let $f(x)=x^2$. Since
$\int_{-\infty}^\infty f(y)P(t\mathrm{e}^{-u},t,x\mathrm{e}^{-\kappa u},\mathrm{d}y) =
P_{t\mathrm{e}^{-u},t} f(x\mathrm{e}^{-\kappa u}) = x^2$
and $A_tf(x) = \frac{2\kappa}{t}x^2$,
the result follows immediately from Equation \eqref{E:QuadVar}.
\end{pf}

Note that $\anglel Z,Z\angler _t = 2\kappa\int_0^t \frac{1}{s} Z_{s}^2
\,\mathrm{d}s$. The predictable quadratic variations of $X$ and $Z$ are given by
the same functional of the process.

%s4.5 #&#
\subsection{Extension to mimicking Brownian martingales}\label{S:Extension}

Now we discuss how we can (and cannot) alter our martingale condition
to mimic some Brownian related processes, including the martingales
associated with the Hermite polynomials and the exponential martingale
of Brownian motion.
%Then we show how Markov projection helps in mimicking continuous
%self-similar martingales if we relax the Markov assumption.
%We finish this section by introducing time-varying self-similar
%processes and extending our construction to mimicking these processes.

Consider the Hermite polynomials $h_n$ which are defined by
\[
\sum_{n\ge0} \frac{u^n}{n!} h_n(x) =
\exp\bigl(ux-u^2/2\bigr)\quad\quad \forall u,x \in \mathbb{R},
\]
equivalently, $h_n(x) = (-1)^n \mathrm{e}^{x^2/2} \frac{\mathrm{d}^n}{\mathrm{d}x^n}(\mathrm{e}^{-x^2/2})$.
Let
\[
H_n(x,t) = t^{n/2}h_n (x/\sqrt{t} )\quad\quad \forall x
\in\mathbb{R}, t>0
\]
and $H_n(x,0)=x^n$.
Then $H_n(B_t,t)$, where $B$ denotes a Brownian motion, is a local
martingale for every $n$ since $H_n(x,t)$ is space--time harmonic,
that is, $\frac{\partial H_n}{\partial t} + \frac{1}{2}\frac{\partial^2
H_n}{\partial x^2} = 0$.

Take $n=2$, the process $H_2(B_t,t) = B_t^2-t$ is Markovian and 1-self-similar,
thus can be mimicked using our mimicking scheme with any subordinator
that satisfies $\psi(1)=1$.

For $n\ge3$, $H_n(B_t,t)$ is $\frac{n}{2}$-self-similar, but it is not
Markovian (see \cite{FanHamKle12}).
So we are not able to mimic this process by a direct application of the
method described above. However, a slight modification of our
construction proves sufficient to achieve our aim.

Let $X_t$ be a mimic of the Brownian motion as in Section~\ref{S:GausContMart} with $k=0$, but without the requirement that $X_t$ be
a martingale. Then we have the following result.

%pr4.6 #&#
\begin{proposition}\label{P:Hermite}
For each $n$, the process $H_n(X_t,t)$ has the same marginal
distributions as $H_n(B_t,t)$ and is a martingale if and only if $\psi
(n/2)=n/2$, or $\mathbb{E} [\mathrm{e}^{-(n/2)\zeta_{\ln(t/s)}} ] =
(s/t)^{n/2}$.
\end{proposition}
\begin{pf}
It is obvious that $H_n(X_t,t)$ and $H_n(B_t,t)$ have the same marginals.
Here we prove the martingale condition.
The transition function of $X$ is
\[
\widetilde{P}(s,t,x,\mathrm{d}y) = \int P\bigl(rt,t,(t/s)^{1/2}r^{1/2}x,\mathrm{d}y
\bigr) G_{s,t}(\mathrm{d}r),\quad\quad s\le t,
\]
where $P$ is the transition function of $B$ and $G_{s,t}$ is the
distribution of $\mathrm{e}^{-\zeta_{\ln(t/s)}}$.
Writing $\mathcal{F}_t$ as the natural filtrations of $H_n(X_t,t)$ and
since $H_n(B_t,t)$ is a martingale, we have
\begin{eqnarray*}
\mathbb{E}\bigl[H_n(X_t,t)|\mathcal{F}_s
\bigr] &=& \mathbb{E} \biggl[\int H_n(y,t)\widetilde{P}(s,t,X_s,\mathrm{d}y)
\big|\mathcal {F}_s \biggr]
\\
&=& \mathbb{E} \biggl[\int H_n \bigl((t/s)^{1/2}r^{1/2}X_s,
rt \bigr) G_{s,t}(\mathrm{d}r) \big|\mathcal{F}_s \biggr].
\end{eqnarray*}
Using that $H_n (ax,t ) = a^nH_n (x,\frac{t}{a^2} )$
for any $a>0$, we then obtain
\[
\mathbb{E}\bigl[H_n(X_t,t)|\mathcal{F}_s
\bigr]  = (t/s)^{n/2} H_n(X_s,s) \mathbb {E}
\bigl[\mathrm{e}^{-(n/2)\zeta_{\ln(t/s)}} \bigr].
\]
Hence, $H_n(X_t,t)$ is a martingale if and only if $\mathbb{E}
[\mathrm{e}^{-(n/2)\zeta_{\ln(t/s)}} ] = (s/t)^{n/2}$.
\end{pf}

Therefore, in order to mimic the process $H_n(B_t,t)$, we can mimic $B_t$,
with the martingale requirement $\psi(\frac{1}{2})=\frac{1}{2}$ changed
to $\psi(\frac{n}{2})=\frac{n}{2}$,
and then apply the function $H_n(x,t)$ to the resulting process.
%A key property in establishing Proposition~\ref{P:Hermite} is the fact
%that $H_n(B_t,t)$ is a martingale.
It is of interest to ask whether the above trick extends to other
space--time harmonic functions. In particular, could this enable us to
mimic the geometric Brownian motion $\exp(B_t-t/2)$. Unfortunately,
this is not the case.
In fact, $H_n(x,t)$, $n\geq1$, are the only analytic functions for
which this trick works.

%pr4.7 #&#
\begin{proposition}
Suppose that $H(x,t) = \sum_{m,n=0}^{\infty} a_{m,n}x^mt^n$ and there
exists $r>1$ such that $\sum_{m,n} |a_{m,n}| r^{m+n} <\infty$, in other
words, $H(x,t)$ is analytic on the set $(-r,r)^2$.
Suppose further that $H$ is space--time harmonic, so that $H(B_t,t)$ is
a martingale. Suppose that $X_t$ mimics $B_t$ with the martingale
requirement replaced with $\psi(k/2)=k/2$ for a positive integer $k$.
Then $H(X_t,t)$ has the same marginals as $H(B_t,t)$ and is a
martingale if and only if $H(x,t)=t^{k/2}h_k(x/\sqrt{t})$ where
$h_k(x)$ is the $k$th Hermite polynomial.
\end{proposition}

\begin{pf}
For $H(X_t,t)$ to be a martingale, we must have for any $x$ and $s\le t$,
\[
H(x,s) = \int H(y,t) \widetilde{P}(s,t,x,\mathrm{d}y) = \int H\bigl((t/s)^{1/2}r^{1/2}x,rt
\bigr) G_{s,t}(\mathrm{d}r),
\]
that is, $\mathbb{E}[H(xu^{1/2}R_u^{1/2},suR_u)] = H(x,s)$ for any $x$,
$s$ and $u\ge1$.
Letting $u=\mathrm{e}^t$ and $V_t = \mathrm{e}^{t-\zeta_t}$, this is equivalent to
$\mathbb{E}[H(xV_t^{1/2}, sV_t)] = H(x,s)$ for all $x$ and $s$, or
\[
\sum_{m,n=0}^\infty a_{m,n}x^ms^n
\mathbb{E} \bigl[V_t^{n+m/2} \bigr] = \sum
_{m,n=0}^\infty a_{m,n} x^ms^n.
\]
Therefore, for all $m$, $n$ and $t<\frac{1}{2}\ln r$, we must have
$a_{m,n}\mathbb{E}[V_t^{n+m/2}] = a_{m,n}$.
Thus, either $a_{m,n}=0$ or $\mathbb{E}[V_t^{n+m/2}]=1$.

Recall that $\mathbb{E}[V_t^\lambda] = \mathbb{E}[\exp(\lambda t-\lambda
\zeta_t)] = \exp(-t(\psi(\lambda)-\lambda))$ and there is at most one
$\lambda$ satisfying $\psi(\lambda)=\lambda$.
Now, choose $\zeta_t$ such that $\mathbb{E}[V_t^{k/2}]=1$ for a $k\in
\mathbb{N}^*$.
Then, for all $(m,n)$ such that $m+2n\neq k$, $a_{m,n}=0$.
Therefore,
\[
H(x,s) = \sum_{n=0}^{\lfloor k/2 \rfloor}
a_{k-2n,n}x^{k-2n}s^n.
\]
Furthermore,
\[
H(cx,s) = c^k \sum_{n=0}^{\lfloor k/2 \rfloor}
a_{k-2n,n} x^{k-2n} \bigl(s/c^2\bigr)^n =
c^k H\bigl(x,s/c^2\bigr).
\]
By Pluci\'nska \cite{Plu98} and Fitzsimmons \cite{Fit01}, $H(x,1)$ is
the $k$th Hermite polynomial.
\end{pf}

%co4.1 #&#
\begin{corollary}
Let $X_t$ be any mimic of $B_t$ in the sense of Section~\ref{S:Scheme}
but without the martingale requirement.
Although the process $\exp(X_t-t/2)$ has the same marginal
distributions as $\exp(B_t-t/2)$, it is not a martingale unless $\zeta
_t =t$, in which case $X=B$.
\end{corollary}
%sA #&#
\begin{appendix}

\section*{Appendix}\label{app}

%leA.1 #&#
\begin{lemma} \label{L:GenTimechange}
Let $(Y_t)_{t\ge0}$ be a Markov process with infinitesimal generator
$A_t$ and $c_t$ be a deterministic, differentiable, increasing function
in $t$ with derivative $c'_t$. Then the time-changed process
$(\widetilde{Y}_t)_{t\ge0} := (Y_{c_t})_{t\ge0}$ is also a Markov
process with infinitesimal generator $\widetilde{A}_t = c'_tA_{c_t}$.
Furthermore, if $f$ is in the domain of $A$, then $f$ is in the domain
of $\widetilde{A}$.
\end{lemma}

\begin{pf}
Let $\mathcal{F}$ be the filtration of $Y$ and $\widetilde{\mathcal
{F}}$ be the filtration of $\widetilde{Y}$. For any bounded measurable
function $g$, we have, for $s\le t$,
\[
\mathbb{E}\bigl[g(\widetilde{Y}_t)|\widetilde{\mathcal{F}}_s
\bigr] = \mathbb {E}\bigl[g(Y_{c_t})|\widetilde{\mathcal{F}}_s
\bigr] = \mathbb {E}\bigl[g(Y_{c_t})|\mathcal{F}_{c_s}\bigr] =
\mathbb{E}\bigl[g(Y_{c_t})|Y_{c_s}\bigr] = \mathbb{E}\bigl[g(
\widetilde{Y}_t)|\widetilde{Y}_s\bigr].
\]
For $t$ where the function $c$ is strictly increasing, the
infinitesimal generator of $Y_{c_t}$ is
\[
\widetilde{A}_tf(x) = \lim_{u\downarrow t}
\frac{\mathbb{E}[f(Y_{c_u}) | Y_{c_t} = x] -
f(x)}{c_u-c_t} \frac{c_u-c_t}{u-t} = A_{c_t}f(x)c'_t.
\]
If $c_u=c_t$ in a small neighbourhood of $t$, then $c'_t=0$ and
\[
\widetilde{A}_tf(x) = \lim_{u\downarrow t}
\frac{\mathbb{E}[f(Y_{c_t})
| Y_{c_t} = x] - f(x)}{u-t} = 0 = A_{c_t}f(x)c'_t.
\]
\upqed\end{pf}

%leA.2 #&#
\begin{lemma} \label{L:GenProduct}
Let $(Y_t)_{t\ge0}$ be a Markov process with infinitesimal generator
$A_t$ and $c_t$ be a deterministic, differentiable function in $t$ with
derivative $c'_t$ and $c_t\neq0$ for any $t$. Then the process
$(\widetilde{Y}_t)_{t\ge0} := (c_tY_t)_{t\ge0}$ is also a Markov
process and has generator
\[
\widetilde{A}_tf(x) = \pi_{ 1/c_t}A_t
\pi_{c_t}f(x) + \frac
{c'_t}{c_t}xf'(x),
\]
where $\pi_c$ is an operator defined by $\pi_cf(x) = f(cx)$.
Furthermore, if $\pi_{c}f$ is in the domain of $A$ for any $c$ and $f$
is differentiable, then $f$ is in the domain of $\widetilde{A}$.
\end{lemma}

\begin{pf}
Let $\mathcal{F}$ be the filtration of $Y$ and $\widetilde{\mathcal
{F}}$ be the filtration of $\widetilde{Y}$.
Let $h_t$ be a function such that $h_t(x)=c_tx$.
Since $h_t$ is one-to-one, $\sigma(h_u(Y_u), u\le s) = \sigma(Y_u,
u\le s)$.
Therefore, for any bounded measurable function $g$, we have
\[
\mathbb{E}\bigl[g(\widetilde{Y}_t)|\widetilde{\mathcal{F}}_s
\bigr] = \mathbb {E}\bigl[g\circ h_t (Y_t)|\widetilde{
\mathcal{F}}_s\bigr] = \mathbb{E}\bigl[g \circ h_t
(Y_t)|\mathcal{F}_s\bigr] = \mathbb{E}\bigl[g\circ
h_t(Y_t)|Y_s\bigr] = \mathbb {E}\bigl[g(
\widetilde{Y}_t)|\widetilde{Y}_s\bigr].
\]
The infinitesimal generator of $\widetilde{Y}$ is
\begin{eqnarray*}
\widetilde{A}_tf(x) &=& \lim_{u\downarrow t}
\biggl(\mathbb{E}\biggl[\pi_{c_u}f(Y_u)\big|Y_t=\frac
{1}{c_t}x\biggr] - \pi_{c_u}f\biggl(\frac{1}{c_t}x\biggr) + f\biggl(\frac{c_u}{c_t}x\biggr) -
f(x)\biggr)\big/(u-t)
\\
&=& A_t\pi_{c_t}f \biggl(\frac{1}{c_t}x \biggr) +
\frac{c'_t}{c_t}xf'(x).
\end{eqnarray*}
\upqed\end{pf}

%leA.3 #&#
\begin{lemma} \label{L:Lamperti}
Suppose $(Z_t)_{t\ge0}$ is a $\kappa$-self-similar Markov process.
Suppose $P(s,t,x,\mathrm{d}y)$ and $L_t$ are, respectively, the transition
function and infinitesimal generator of $Z$.
Let $\widehat{Z}_t = \mathrm{e}^{-t\kappa}Z_{\mathrm{e}^t}$. Then $(\widehat{Z}_t)_{t\in
\mathbb{R}}$ is a time-homogeneous Markov process with transition semigroup
\[
\widehat{P}_tf(x) = \int f(y) P\bigl(\mathrm{e}^{-t},1,x\mathrm{e}^{-t\kappa},\mathrm{d}y
\bigr)
\]
and infinitesimal generator
\[
\widehat{L}f(x) = L_1 f(x) - {\kappa}xf'(x).
\]
Furthermore, if $f$ is in the domain of $L$ and differentiable, then it
is in the domain of $\widehat{L}$.
\end{lemma}

\begin{pf}
By the scaling property of $P$, we have
\[
\mathbb{P} (\widehat{Z}_t\in \mathrm{d}y | \widehat{Z}_s=x ) = P
\bigl(\mathrm{e}^s, \mathrm{e}^t, x\mathrm{e}^{s\kappa}, \mathrm{e}^{t\kappa}\,\mathrm{d}y
\bigr) = P\bigl(\mathrm{e}^{-(t-s)}, 1, x\mathrm{e}^{-(t-s)\kappa}, \mathrm{d}y\bigr).
\]
It follows that $\widehat{Z}$ is time-homogeneous and $\widehat
{P}_tf(x) = \int f(y) P(\mathrm{e}^{-t},1,x\mathrm{e}^{-t\kappa},\mathrm{d}y)$.
The generator of $\widehat{Z}$ can be obtained by applying Lemma~\ref
{L:GenTimechange} and Lemma~\ref{L:GenProduct},
and seeing that $\pi_{\mathrm{e}^{t\kappa}}\mathrm{e}^tL_{\mathrm{e}^t}\pi_{\mathrm{e}^{-t\kappa}}=L_1$
from Equation \eqref{E:SSGen} with $c=\mathrm{e}^t$ and $s=1$.

Note that for all $c>0$, $\pi_c f$ is in the domain of $L$ by the
scaling property of $Z$.
Thus, writing $\check{L}$ as the generator of $\check{Z}_t := Z_{\mathrm{e}^t}$,
$\pi_c f$ is in the domain of $\check{L}$ by Lemma~\ref{L:GenTimechange}.
\end{pf}

%leA.4 #&#
\begin{lemma} \label{L:GenBochner}
Suppose $(\chi_t)_{t\ge0}$ is a time-homogeneous Markov process with
semigroup $P_t$ and generator $L$, and $(\eta_t)_{t\ge0}$ is a
subordinator independent of $\chi$ with drift $\beta$ and L\'evy
measure $\nu$. Set $Y_t = \chi_{\eta_t}$. Then the process $(Y_t)_{t\ge
0}$ is a time-homogeneous Markov process with generator $A$ where
\[
Af(x) = \beta Lf(x) + \int_{(0,\infty)} \bigl(P_uf(x)-f(x)
\bigr) \nu(\mathrm{d}u).
\]
Furthermore, if $f$ is in the domain of $L$, then it is in the domain
of $A$.

If $\eta$ has zero drift, then $Y$ is a pure jump process.
\end{lemma}

\begin{pf}
See Sato \cite{Sat99}, Theorem~32.1.
\end{pf}

%\begin{appendix}
%\section{}
\end{appendix}

\section*{Acknowledgements}
The authors are grateful to two referees and an associate editor for
their careful reading of an earlier version of the
paper and a number of suggestions and improvements.

This research was supported by the
Australian Research Council Grant DP0988483 and the first author was
the recipient of a Victorian International Research Scholarship.
%\begin{supplement}%[id=suppA]
%\slink[doi]{10.3150/00-BEJXXXXSUPP} %[doi,text={...}] - jei reikia
%suskaldyti doi
%\sfilename{BEJ000\_supp.pdf}
%\end{supplement}

% imsref loaded by arune.pranskunaite, 2014-03-06 13:04:51
%
% imsref loaded by arune.pranskunaite, 2014-03-07 14:56:40

\printhistory

\end{document}